\documentclass[a4paper,11pt]{article}

\usepackage{geometry,enumerate,amsmath,amssymb,latexsym}
\usepackage[all]{xy}

\newtheorem{example}{Example}[section]

\newtheorem{defn}[example]{Definition}

\newtheorem{prop}[example]{Proposition}

\newtheorem{thm}[example]{Theorem}

\newtheorem{rem}[example]{Remark}
\newtheorem{remark}[example]{Remark}
\newtheorem{cor}[example]{Corollary}

\newenvironment{pf}{{\bf Proof:}}{\hfill $\Box$

\mbox{}}

\def\geq{\geqslant}

\def\leq{\leqslant}

\def\ge{\geqslant}

\def\io{^{-1}}
\def\A{\alpha}
\begin{document}

\title{\Large \bf  LIE LOCAL SUBGROUPOIDS AND THEIR \\
HOLONOMY AND MONODROMY LIE GROUPOIDS\thanks{KEYWORDS: local
equivalence relation, local subgroupoid, holonomy groupoid,
monodromy groupoid, monodromy principle: 1991 AMS Classification:
58H05,22A22,18F20} }

\small{ \author  {Ronald Brown  \\ School of Informatics \\
Mathematics Division
\\ University of Wales  \\ Bangor, Gwynedd \\ LL57 1UT, U.K.
\\r.brown@bangor.ac.uk \\ \and
\.{I}lhan \.{I}\c{c}en    \\  University of  \.{I}n\"{o}n\"{u} \\
Faculty of Science and Art
\\ Department of
Mathematics
\\ Malatya/ Turkey \\ iicen@inonu.edu.tr}}
\maketitle
\begin{abstract}
The notion of local equivalence relation on a topological space is
generalised to that of local subgroupoid. The main result is the
construction  of the holonomy and monodromy groupoids of certain
Lie local subgroupoids, and the formulation of a monodromy
principle on the extendibility of local Lie morphisms.
\end{abstract}

\section*{Introduction}

It has long been recognised that the notion of Lie group is
inadequate to express the local-to-global ideas inherent in the
investigations of Sophus Lie, and various extensions have been
developed, particularly the notion of Lie groupoid, in the hands
of Ehresmann, Pradines, and others.

Another set of local descriptions have been given in the notion of
foliation (due to Ehresmann) and also in the notion of local
equivalence relation (due to Grothendieck and Verdier).

Pradines in \cite{Pr1} also introduced the notion of what he
called `morceau d'un groupo\"{\i}de de Lie' and which we have
preferred to call `locally Lie groupoid' in \cite{Br-Mu1}. This is
a groupoid $G$ with a subset $W$ of $G$ containing the identities
of $G$ and with a manifold structure on $W$ making the structure
maps `as smooth as possible'. It is a classical result that in the
case $G$ is a group the manifold structure can be transported
around $G$ to make $G$ a Lie group. This is false in general for
groupoids, and this in fact gives rise to the {\it holonomy
groupoid} for certain such $(G,W)$.

In \cite{Br-Mu2} it is shown that a foliation on a paracompact
manifold gives rise to a locally Lie groupoid. It is part of the
theory of Lie groupoids that a Lie algebroid gives rise, under
certain conditions, to a locally Lie groupoid. Thus a locally Lie
groupoid is one of the ways of giving a  useful expression of
local-to-global structures.

The notion of {\em local equivalence relation}  was introduced by
Grothendieck and Verdier \cite{Gr-Ve} in a series of exercises
presented as open problems concerning the construction of a
certain kind of topos. It  was investigated further by Rosenthal
\cite{Ro1,Ro2} and more recently by Kock and Moerdijk
\cite{Ko-Mo1,Ko-Mo2}. A local equivalence relation is a global
section of the  sheaf ${\mathcal E}$  defined by the presheaf $E$
where $E(U)$ is the set of all equivalence relations on the open
subsets $U$ of $X$,  and $E_{UV}$ is the restriction map from
$E(U)$ to $E(V)$ for $V\subseteq U$. The main aims of the papers
\cite{Gr-Ve,Ko-Mo1,Ko-Mo2,Ro1,Ro2} are towards the connections
with sheaf theory and topos theory. Any foliation gives rise to a
local equivalence relation, defined by the path components of
local intersections of small open sets with the leaves.

An equivalence relation on a set $U$ is just a wide subgroupoid of
the indiscrete groupoid $U \times U$ on $U$. Thus it is natural to
consider the generalisation which replaces the indiscrete groupoid
on the topological space $X$ by any groupoid $Q$ on $X$. So we
define a {\em local subgroupoid} of the  groupoid $Q$ to be a
global section of the sheaf ${\mathcal L}$ associated to the
presheaf $L_Q$  where $L(U)$ is the set  of all wide subgroupoids
of $Q|U$ and $L_{UV}$ is the restriction map from $L(U)$ to $L(V)$
for $V\subseteq U$.

Our aim is towards local-to-global principles and in particular the
monodromy principle, which in our terms is formulated as the
globalisation of local morphisms (compare \cite{Ch,Pr1,Br-Mu1}).
Our first formulation is for the case $Q$ has no topology, and this
gives our `weak monodromy principle' (Theorem \ref{wmp}).

In the case $Q$ is a Lie groupoid we expect to deal with Lie local
subgroupoids $s$ and the globalisation of local smooth morphisms
to a smooth morphism $Mon(s)\to K$ on a `monodromy Lie groupoid'
$Mon(s)$ of $s$. The construction of the Lie structure on $Mon(s)$
requires extra conditions on $s$ and its main steps are:

$\bullet$ \   the construction of  a locally Lie groupoid from $s$
and a strictly regular atlas for $s$,

$\bullet$ \  applying the construction of the holonomy Lie
groupoid of the locally Lie groupoid, as in \cite{Ao-Br,Br-Mu1},

$\bullet$ \  the further construction of the monodromy Lie groupoid,
as in \cite{Br-Mu2}.

For strictly regular atlases ${\cal U}_s= \{(U_i, H_i):i\in I\}$
for $s$ this leads to a morphism of Lie groupoids
\[   \zeta : Mon(s, {\cal U}_s)\to Hol(s, {\cal U}_s) \]
each of which contains the $H_i, i\in I$, as   Lie subgroupoids,
and which are in a certain sense maximal and minimal respectively
for this property.  This morphism $\zeta$ is \'etale on stars.
Further, a smooth local morphism $\{f_i:H_i\to K,i\in I\}$ to a
Lie groupoid $K$ extends uniquely to a smooth morphism $Mon(s,
{\cal U}_s)\to K$. This is our strong monodromy principle (Theorem
\ref{smp}).

It should be noticed that this route to a monodromy Lie groupoid
is different from that commonly taken in the theory of foliations.
For a foliation $\mathcal{F}$ it is possible to define the
monodromy groupoid as the union of the fundamental groupoids of
the leaves, and then to take the holonomy groupoid as a quotient
groupoid of this, identifying classes of paths which induce the
same holonomy.

However there seem to be strong advantages in seeing these
holonomy and monodromy groupoids as special cases of much more
general constructions, in which the distinct universal properties
become clear. In particular, this gives a link between the
monodromy groupoid and the important monodromy principle, of
extendability of local morphisms. In the Lie case, this requires
moving away from the \'etale groupoids which is the main emphasis
in  \cite{Ko-Mo1,Ko-Mo2}.

We hope to investigate elsewhere the relation of these ideas to
questions on transformation groups.

\section{Local Subgroupoids}

Consider a groupoid $Q$ on a set $X$ of objects, and suppose also
$X$ has a topology. For any open subset $U$ of $ X$ we write $Q|U$
for the full subgroupoid of $Q$ on the object set $U$. Let
$L_Q(U)$ denote the set of all wide subgroupoids of $ Q|U$. For
$V\subseteq U$, there is a restriction map $L_{UV}\colon L_Q(U)
\to L_Q(V)$ sending $H$ in $L_Q(U)$ to $H|V$. This gives $L_Q$ the
structure of presheaf on $X$.

We first  interpret in our case the usual  construction of the
sheaf $p_Q:{\cal L}_Q \to X$ constructed from the presheaf $L_Q$.

For $x\in X$, the stalk ${p_Q}^{-1}(x)$ of ${\cal L}_Q$ has
elements the germs $[U, H_U]_x$ where $U$ is open in $X$, $x\in
U$, $H_U$ is a wide subgroupoid of $Q|U$, and the equivalence
relation $\sim_x$ yielding the germs at $x$ is that $H_U\sim_x
K_V$, where $K_V$ is wide subgroupoid of $Q|V$, if and only if
there is a neighbourhood $W$ of $x$ such that $W\subseteq U\cap V$
and $H_U|W =K_V|W$.

\begin{defn}{\rm
A {\it local subgroupoid} of $Q$ on the topological space $X$ is a
global section of the sheaf $p_Q:{\cal L}_Q\to X$ associated to
the presheaf $L_Q$.}
\end{defn}

An {\it atlas } ${\cal U}_s = \{(U_i, H_i):i\in I\}$ for a local
subgroupoid $s$ of $Q$ consists of an open cover ${\cal U}=
\{U_i:i\in I\}$ of $X$, and for each $i\in I$ a wide subgroupoid
$H_i$ of $Q|U_i$ such that for all $x\in X$, $i\in I$, if $x\in
U_i$ then $s(x)=[U_i, H_i]_x$.

Two standard  examples of $Q$ are   $Q=X$, $Q=X\times X$.  In the
first case, $L_X$ is a sheaf and ${\cal L}_X\to X$ is a bijection.
In the case $Q$ is the indiscrete groupoid $X\times X$ with
multiplication $(x, y)(y, z) = (x, z)$, $x, y, z\in X$, the local
subgroupoids of $Q$ are the local equivalence relations on $X$, as
mentioned in the Introduction. It is known that $L_{X\times X}$ is
in general not a sheaf \cite{Ro1}.

In the following, we show that many of the basic results obtained
by Rosenthal in \cite{Ro1,Ro2} extend conveniently to the local
subgroupoid case.

The set $L_Q(X)$ of wide subgroupoids of $Q$ is a poset under
inclusion. We write $\leq$ for this partial order.

Let ${\bf Loc}(Q)$ be the set of local subgroupoids of $Q$. We
define a partial order $\leq$ on ${\bf Loc}(Q)$ as follows.

Let $x\in X$. We define a partial order on the stalks
${p_Q}^{-1}(x)={{\cal L}^Q}_x$ by $[U', H']_x\leq [U, H]_x$ if
there is an  open neighbourhood $W$ of $x$ such that $W\subseteq
U\cap U'$ and $H'|W$ is a subgroupoid of $H|W$. Clearly this
partial order is well defined. It induces a partial order on ${\bf
Loc}(Q)$ by $s\leq t$ if and only if $s(x)\leq t(x)$ for all $x\in
X$.

We now fix a groupoid  $Q$  on $X$, so that $L_Q(X)$ is the set of
wide subgroupoids of $Q$, with its inclusion partial order, which
we shall write $\leq$.

We define poset morphisms
\[ loc_Q : L_Q(X)\to {\bf Loc}(Q)    \ \ \mbox{and}\
\ \ \ \  glob_Q:{\bf Loc}(Q)\to L_Q(X) \] as follows. We
abbreviate $loc_Q$, $glob_Q$ to $loc, glob$.
\begin{defn}{\rm
If $H$ is a wide subgroupoid of the groupoid $Q$ on $X$, then
$loc(H)$ is the local subgroupoid defined by
\[             loc(H)(x) =[X, H]_x.    \]
Let $s$ be a local subgroupoid of  $Q$. Then $glob(s)$ is the wide
subgroupoid of $Q$ which is the intersection  of all wide
subgroupoids $H$ of $Q$ such that $s\leq loc(H)$.} \end{defn}

We think of $glob(s)$ as an approximation to $s$ by a global
subgroupoid.

\begin{prop}\label{sim}

i) $loc$ and $glob$ are morphisms of posets.

ii) For any wide subgroupoid $H$ of $Q$, $glob(loc(H))\leq H$.
\hfill $\Box$
\end{prop}
The proofs are clear.

However, $s\leq loc(glob(s))$ need not hold. Examples of this are
given in Rosenthal's paper \cite{Ro1} for the case of local
equivalence relations.

Here is an alternative description of $glob$. Let ${\cal U}_s =
\{(U_i, H_i):i\in I\}$ be an atlas for the local subgroupoid $s$.
We define $glob({\cal U}_s)$ to be the subgroupoid of $Q$
generated by all the $H_i, i \in I$.

An atlas ${\cal V}_{s}=\{(V_j, {s}_j): j\in J\}$ for $s$ is
said to refine ${\cal U}_s$ if for each index $j\in J$
there exists an index $i(j)\in I$ such that $V_j\subseteq U_{i(j)}$
and $s_{i(j)}|V_j = s_j$.

\begin{prop} \label{refin2}
Let  $s$  be a local subgroupoid of $Q$ given by the atlas ${\cal
U}_s = \{(U_i, H_i):i\in I\}$. Then $glob(s)$ is the intersection
of the subgroupoids $glob({\cal V}_s) $ of $Q$ for all refinements
${\cal V}_s$ of ${\cal U}_s$.
\end{prop}
\begin{pf}
Let $K$ be the intersection given in the proposition.

Let $S$ be a subgroupoid of $Q$ on $X$ such that $s\leq loc(S)$.
Then for all $x \in X$ there is a neighbourhood $V$ of $x$ and
$i_x \in I$ such that $x \in U_{i_x}$ and $ H_{i_x}|V_x \cap
U_{i_x} \leq S$. Then ${\cal W}= \{(V_x \cap U_{i_x},H_{i_x}|V_x
\cap U_{i_x}) : x \in X\} $ refines ${\cal U}_s$ and $glob({\cal
W}) \leq S$. Hence $K \leq S$, and so $ K \leq glob(s) $.

Conversely, let ${\cal V}_s = \{ (V_j,H'_j) : j \in J \}$ be an
atlas for $s$ which refines ${\cal U}_s$. Then for each $ j \in J$
there is an $i(j) \in I$ such that $ V_j \subseteq
U_{i(j)},H'_j=H_{i(j)}|V_j$. Then $s \leq loc(glob({\cal V}_s))$. Hence
$glob(s) \leq glob({\cal V}_s)$ and so $glob(s) \leq K$.
\end{pf}

We need the next definition in the following sections.

\begin{defn}\label{def1}{\em
Let $s$ be a local subgroupoid of the groupoid $Q$ on $X$. An
atlas ${\cal U}_s$ for $s$ is called {\it  globally adapted} if $
glob(s)= glob({\cal U}_s)$.}
\end{defn}

\begin{rem}{\em This is a variation on the notion of an
$r$-adaptable family defined by Rosenthal in \cite[Definition
4.4]{Ro2} for the case of a local equivalence relation $r$. He
also imposes a connectivity condition on the local equivalence
classes.     }   \end{rem}

\section{The weak monodromy principle for local subgroupoids} \label{wmonod}

Let  $s$  be a  local subgroupoid of $Q$ which is given by an
atlas ${\cal U}_s = \{(U_i, H_i): i\in I\}$, and let $H=glob(s) $,
$W({\cal U}_s)=\bigcup_{i\in I}H_i$. Then $ W({\cal U}_s)\subseteq
H$.

The set $W({\cal U}_s)$  inherits  a {\it pregroupoid}  structure
from  the groupoid $H$. That is, the source and target maps
$\alpha, \beta$ restrict to maps on $W({\cal U}_s)$, and if $u,
v\in W({\cal U}_s)$ and $\beta u = \alpha v$, then the composition
$uv$ of $u,v$ in $H$ may or may not belong to $W({\cal U}_s)$. We
now follow the method of Brown and Mucuk in \cite{Br-Mu1}, which
generalises work for groups in Douady and Lazard \cite{Do-La}.

There is a standard construction $M(W({\cal U}_s))$  associating
to the pregroupoid  $W({\cal U}_s)$  a morphism
$\tilde{\imath}\colon W({\cal U}_s)\rightarrow M(W({\cal U}_s))$
to a groupoid $M(W({\cal U}_s))$ and which  is  universal for
pregroupoid morphisms to a groupoid.  First, form the free
groupoid $F(W({\cal U}_s))$  on the graph  $W({\cal U}_s)$, and
denote the inclusion $W({\cal U}_s)\rightarrow F(W({\cal U}_s))$
by  $u\mapsto [u] $. Let   $N$  be  the normal  subgroupoid
(Higgins~\cite{Hi}, Brown~\cite{Br1})  of $F(W({\cal U}_s))$
generated by the elements $[vu]^{-1}[v][u]$  for all $u,v\in
W({\cal U}_s)$   such  that $vu$   is defined and belongs to
$W({\cal U}_s) $. Then $M(W({\cal U}_s))$  is defined to be  the
quotient  groupoid (loc. cit.)   $F(W({\cal U}_s))/N $.  The
composition   $W({\cal U}_s)\rightarrow F(W({\cal
U}_s))\rightarrow M(W({\cal U}_s))$ is written $\tilde{\imath} $,
and is the required universal morphism.

There is a unique morphism of groupoids  $p\colon M(W({\cal
U}_s))\rightarrow glob(s)$  such  that $p\tilde{\imath}$  is the
inclusion $i\colon W({\cal U}_s)\rightarrow glob(s) $. It follows
that  $\tilde{\imath}$  is injective.  Clearly,  $p$ is surjective
if and only if  the atlas for $s$ is globally adapted.  In this
case, we call  $M(W({\cal U}_s))$  the {\it monodromy groupoid} of
$W({\cal U}_s)$ and write it $Mon(s, {\cal U}_s)$.

\begin{defn}
{\em The local subgroupoid $s$ is called {\it simply connected }
if it has a globally  adapted atlas ${\cal U}_s$ such that the
morphism $p: Mon(s,{\cal U}_s) \to glob(s)$ is an isomorphism.}
\end{defn}

We now relate $Mon(s,{\cal U}_s)$ to the extendability of local
morphisms to a groupoid $K$.

Let $K$ be a groupoid with object space $X$.

\begin{defn}{\em
A {\it local morphism}
$\xy(0,0)*+{f:{\cal U}_s}\ar @{-->}(17,0)*+{K}
\endxy$
consists of a globally adapted atlas
${\cal U}_s =\{(U_i, H_j):i\in I\}$
for $s$ and a
family of morphisms $f_i: H_i\to K, i\in I$ over the inclusion
$U_i\to X$ such that for all $i,j\in I$,
\[  f_i|(H_i\cap H_j) = f_j|(H_i\cap H_j), \]
and  the resulting function $f':W({\cal U}_s)\to K$ is a
pregroupoid morphism.}
\end{defn}

\begin{thm}\label{wmp}(Weak Monodromy Principle) A local morphism
$\xy(0,0)*+{f:{\cal U}_s}\ar @{-->}(17,0)*+{K} \endxy$ defines
uniquely a groupoid morphism $M(f):Mon(s,{\cal U}_s)\to K$ over
the identity on objects such that $M(f)|H_i = f_i, i\in I$.
Further, if $s$ is simply connected, then the $(f_i)$ determine a
groupoid morphism $glob(s)\to K$.
\end{thm}
\begin{pf}
The proof is direct from the definitions. A local morphism $f$
defines a pregroupoid morphism $f': W({\cal U}_s)\to K$ which
therefore defines $M(f):Mon(s,{\cal U}_s)\to K$ by the universal
property of $W({\cal U}_s)\to Mon(s,{\cal U}_s)$.
\end{pf}

In the next section, we will show how to extend this result to the
Lie case. This involves discussing the construction of a topology
on $Mon(s,{\cal U}_s)$ under the given conditions. For this we
follow the procedure of Brown-Mucuk in \cite{Br-Mu1} in using the
construction and properties of the holonomy groupoid of a locally
Lie groupoid. This procedure is in essence due to Pradines, and
was announced without detail in \cite{Pr1}. As explained in the
preliminary preprint \cite{Br2} these details were communicated by
Pradines to Brown in the 1980s.

\section{Local Lie  subgroupoids,   holonomy  and monodromy}
\label{loclie}

The aim of this section is to give sufficient conditions on local
subgroupoid $s$ of $G$ for the monodromy groupoid of $s$ to admit
the structure of a Lie groupoid, so that the globalisation
$f:Mon(s,{\cal U}_s)\to K$ of a local smooth morphism
$\xy(0,0)*+{f_i:H_i}\ar @{-->}(17,0)*+{K} \endxy$, $i\in I$, is
itself smooth. As explained in the Introduction, our method
follows \cite{Br-Mu1} in first constructing  a locally Lie
groupoid $(glob(s), W({\cal U}_s))$;   the holonomy Lie groupoid
of this locally Lie groupoid comes with a morphism of groupoids
$\psi :Hol(glob(s), W({\cal U}_s))\to glob(s)$ which is a minimal
smooth overgroupoid of $glob(s)$ containing $W({\cal U}_s)$ as an
open subspace. From this holonomy Lie groupoid we construct the
Lie structure  on the monodromy groupoid. We begin therefore by
recalling  the holonomy groupoid construction.

We consider ${\cal C}^r$-manifolds for $r\geq-1$. Here a
${\cal C}^{-1}$-manifold is simply a topological space and for
$r=-1$, a smooth map is simply a continuous map. Thus the Lie
groupoids in the ${\cal C}^{-1}$ case will simply be the
topological groupoids. For $r=0$, a ${\cal C}^0$-manifold is as
usual a topological manifold, and a smooth map is just a
continuous map. For $r\geq 1$, $r=\infty, \omega$ the
definition of ${\cal C}^r$-manifold and smooth map are as
usual. We now fix $r\geq -1$.

One of the key differences between the cases $r = -1$ or $0$ and
$r \ge 1$ is that for $r \ge 1$, the pullback of ${\cal C}^r $
maps need not be a smooth submanifold of the product, and so
differentiability of maps on the pullback cannot always be
defined. We therefore adopt the following definition of Lie
groupoid. Mackenzie~\cite{Ma} (pp. 84-86) discusses the utility
of various definitions of differentiable groupoid.

Recall that if  $G$  is a groupoid then the difference map on $G$
is $\delta :  G \times _{\alpha} G \rightarrow G , (g,h)
\mapstochar \rightarrow g^{-1}h $.

A {\em Lie groupoid} is a topological  groupoid $G$ such that
\begin{enumerate}[(i)]
\item  the space of arrows is a smooth  manifold,  and
the space of objects is a smooth  submanifold of $G$,
\item the  source and target maps $\alpha , \beta $, are smooth
maps and are submersions,
\item  the domain $G \times
_{\alpha } G$ of the difference map $\delta$ is a  smooth
submanifold of $G \times G$,
\item  the difference map $\delta $ is
a  smooth  map.
\end{enumerate}
The term {\em locally Lie groupoid} $(G, W)$ is defined later.

The following definition is due to Ehresmann~\cite{Eh1}.

\begin{defn} {\em Let  $G$  be a groupoid and let  $X = O_{G}$
be a smooth  manifold. An {\em admissible local section} of $G$ is
a function $\sigma : U \rightarrow G$  from an open set in $X$
such that \begin{enumerate}[(i)] \item $\alpha \sigma(x) = x$ for
all $x \in U$; \item  $\beta \sigma(U)$ is open in $X$, and
\item  $\beta \sigma$ maps $U$  diffeomorphically to $\beta \sigma(U)$.
\end{enumerate}}  \label{locsecdef}

\end{defn}

Let $W$ be a subset of $G$ and let $W$ have the structure of a
smooth manifold such that $X$ is a submanifold. We say that
$(\alpha ,\beta ,W)$  is {\em  locally sectionable} if for each $w
\in W$ there is an admissible  local section $\sigma : U
\rightarrow G$ of $G$ such that (i) $\sigma\alpha (w) = w$,  (ii)
$\sigma(U) \subseteq W$ and (iii) $\sigma$ is smooth as a function
from $U$  to $W$. Such a $\sigma$ is called a {\it smooth
admissible local section.}

The following definition is due to Pradines~\cite{Pr1} under the
name  ``{\it morceau de groupoide diff\'erentiables}''.
\begin{defn} {\em A {\it locally Lie  groupoid} is a  pair $(G,W)$
consisting of a groupoid $G$ and a smooth manifold  $W$ such that:
\begin{enumerate}[G1)]
\item $O_G \subseteq W \subseteq G $;
\item  $W = W^{-1} $;
\item  the set $W(\delta ) = (W \times _{\alpha} W) \cap \delta ^{-1}(W)$
 is open in $W \times _{\alpha} W$ and
the restriction of $\delta $ to  $W(\delta )$ is smooth;
\item the restrictions to $W$ of the source and
target maps $\alpha $  and $\beta $ are smooth and the triple
$(\alpha ,\beta ,W)$ is  locally sectionable;
\item $ W$ generates $G$ as a groupoid.
\end{enumerate}}
\end{defn}

Note that in this definition, $G$ is a groupoid but does not need
to have a topology. The locally Lie groupoid $(G,W)$ is said  to
be {\it extendable} if there can be found a topology on $G$ making
it a Lie groupoid and for which $W$ is an open submanifold.  In
general, $(G, W)$ is not extendable, but there is a holonomy
groupoid $Hol(G, W)$ and a morphism $\psi: Hol(G, W)\to  G$ such
that $Hol(G, W)$ admits the structure of Lie groupoid and is the
``minimal'' such overgroupoid of $G$. The construction is given in
detail in \cite{Ao-Br} and is outlined below.

\begin{defn}{\em
A {\it  Lie local  subgroupoid} $s$ of a Lie groupoid $Q$ is a
local subgroupoid $s$ given by an atlas ${\cal U}_s= \{(U_i, H_i):
i\in I\}$ such that for $i\in I$ each $H_i$ is a  Lie subgroupoid
of $Q$.}
\end{defn}
We know from  examples for foliations and hence for local
equivalence relations that $glob(s)$  need not  be  a Lie
subgroupoid of $Q$ \ \cite{Br-Mu2}. Our aim is to define a
holonomy groupoid $Hol(s,{\cal U}_s)$ which is a Lie groupoid.

We now adapt some definitions from \cite{Ro2}.

\begin{defn}{\em
An atlas ${\cal U}_s = \{(U_i, H_i) : {i\in I}\}$ for a Lie local
subgroupoid $s$ of $Q$ is said to be {\it regular} if the groupoid
$(\alpha_i, \beta_i, H_i)$  is locally sectionable for all ${i\in
I}$. A Lie local subgroupoid $s$ is {\em regular} if it has a
regular atlas.}
\end{defn}

\begin{defn}{\em
An atlas ${\cal U}_s=\{(U_i, H_i):{i\in I}\}$ for  a Lie local
subgroupoid $s$ is  said to be  {\it strictly regular } if
\begin{enumerate}[i)]
\item ${\cal U}_s$ is globally adapted to $s$,
\item ${\cal U}_s$ is regular,
\item $W({\cal U}_s)$ has with its topology  as a subset of $Q$  the
structure of smooth submanifold containing each $H_i$, $i\in I$,
as an open submanifold of  $W({\cal U}_s)$ and such that $W({\cal
U}_s)(\delta)$ is open in $W({\cal U}_s)\times_{\alpha} W({\cal
U}_s)$.
\end{enumerate}
A  Lie local subgroupoid $s$ is {\em strictly regular} if it has a
strictly regular  atlas.}
\end{defn}

\begin{rem}{\em The main result of  \cite{Br-Mu2} is that
the local equivalence relation defined by a foliation on
a paracompact manifold has a strictly regular atlas. }
\end{rem}

The following is  a key construction of a locally Lie groupoid from
a strictly regular  Lie local subgroupoid.

\begin{thm}
Let $Q$ be a Lie groupoid on $X$ and let ${\cal U}_s =\{(H_i,
U_i):{i\in I}\}$ be a  strictly regular atlas for the Lie  local
subgroupoid  $s$  of  $Q$. Let
\[           G = glob(s) , \ \  \ \   \  \  \  \    W({\cal U}_s)
=\bigcup_{i\in I} H_i.    \] Then $(G, W({\cal U}_s))$ admits the
structure of  a locally Lie groupoid.
\end{thm}
\begin{pf}

\noindent (G1)\ By the definition of $G$ and $W({\cal U}_s)$, 
clearly  $X\subseteq W({\cal U}_s)\subseteq H$.

\noindent (G2) In fact, $W({\cal U}_s) = W({\cal U}_s)^{-1}$. Let
$g\in W({\cal U}_s)$. Then there is an index $i\in I$ such that
$g\in H_i$. Since $H_i$ is a groupoid on $U_i$, \ $g^{-1}\in H_i$.
So $W({\cal U}_s) = W({\cal U}_s)^{-1}$.

\noindent (G3)  Since $s$ is strictly regular, by definition,
$W({\cal U}_s(\delta))$ is open in
$W({\cal U}_s)\times_{\delta} W({\cal U}_s)$.

We now prove the restriction of $\delta$ to $W({\cal U}_s)(\delta)$ is
smooth.

For each $i\in I$, $H_i$ is a Lie groupoid on $U_i$ and
so the difference map
\[               \delta_i : H_i\times_{\alpha} H_i \rightarrow H_i             \]
is smooth. Because  $H_i\subseteq W({\cal U}_s), \ i\in I$,
using the smoothness of the inclusion
map $i_{H_i}: H_i\rightarrow W({\cal U}_s)$, we get a smooth map
\[    i_{H_i}\times i_{H_i} : H_i\times_{\alpha }H_i\rightarrow
W({\cal U}_s)\times_{\alpha} W({\cal U}_s).    \]
The restriction of $W({\cal U}_s)({\delta})$ is also smooth, that is,
\[         i_{H_i}\times i_{H_i} : H_i\times_\alpha H_i\rightarrow
W({\cal U}_s)(\delta)    \] is smooth. Then the following diagram
is commutative: $$ \xymatrix{ H_i\times_{\alpha} H_i   \ar [r]
^-{\delta}\ar [d]_{i_{H_i}\times i_{H_i}}
 & H_i   \ar [d]^{i_{H_i}}    \\
 W({\cal U}_s)({\delta}) \ar [r] _-{\delta}      &W({\cal U}_s)}
$$
This verifies (G3), since $H_i$ is open in $W({\mathcal U}_s)$ and hence
$H_i \times _{\alpha} H_i$ is open in $W({\mathcal U}_s)(\delta)$.

\noindent (G4) We define source and target maps $\alpha_{W({\cal U}_s)}$
and
$\beta_{W({\cal U}_s)}$ respectively as
follows: if $g\in W({\cal U}_s)$ there exist $i\in I$
such that $g\in H_i$  and we let
\[      \alpha_{W({\cal U}_s)}(g) = \alpha_i(g), \  \   \beta_{W({\cal U}_s)}(g) =
\beta_i(g)       \]
Clearly $\alpha_{W({\cal U}_s)}$ and $\beta_{W({\cal U}_s)}$ are smooth.
Since ${\cal U}_s = \{(U_i, H_i):{i\in I}\}$
is   strictly regular,  $(\alpha_i, \beta_i, H_i)_{i\in I}$
is locally sectionable for all $i\in I$.
Hence $(\alpha_{W({\cal U}_s)}, \beta_{W({\cal U}_s)}, W({\cal U}_s))$
is locally sectionable.

\noindent (G5)  Since the atlas ${\cal U}_s$ is globally adapted
to  $s$, then $G=glob(s)$ is generated by the  $\{H_i\}$, ${i\in
I}$, and so is also generated by  $W({\cal U}_s)$.

Hence $(glob(s), W({\cal U}_s))$ is a locally Lie groupoid.
\end{pf}

There is a main globalisation theorem for a locally topological
groupoid due to Aof-Brown \cite{Ao-Br}, and a Lie version of this
is stated by Brown-Mucuk \cite{Br-Mu1}; it shows how a locally Lie
groupoid gives rise to its holonomy groupoid, which is a Lie
groupoid satisfying a universal property. This theorem gives a
full statement and proof of a part of Th\'{e}or\`eme 1 of
\cite{Pr1}. We can give immediately the generalisation to Lie
local subgroupoids.

\begin{thm} { (Globalisability Theorem)} \  Let  $s$
be  a  Lie local subgroupoid of  a Lie groupoid  $Q$,  and suppose
given a  strictly regular atlas ${\cal U}_s = \{(U_i, H_i): i \in
I\}$ for $s$.  Let $(glob(s), W({\cal U}_s))$  be the associated
locally Lie groupoid.  Then there is a Lie groupoid
$Hol=Hol(s,{\cal U}_s)$,  a morphism $\psi : Hol \rightarrow
glob(s)$   of groupoids and an embedding $i_s: W({\cal
U}_s)\rightarrow Hol$ {\it of} $W({\cal U}_s)$   to an open
neighborhood of $O_{Hol}$ in   $Hol$  such that the following
conditions are satisfied:

\noindent i) $\psi $  is the identity on object, $\psi i_s =
id_{W({\cal U}_s) },  \psi ^{-1}(H_i)$  is open in $Hol$, and the
restriction  $\psi _{H_i} : \psi^{-1}(H_i)  \rightarrow H_i$  of
$\psi $  is smooth;

\noindent  ii) Suppose $A$ is a Lie  groupoid on $X=Ob(Q)$ and
$\xi
: A \rightarrow glob(s)$ is a morphism of groupoids such that:
\begin{enumerate}[a)]
\item $\xi $ is the identity on objects;
\item for all $i$ the restriction $\xi _{H_i} : \xi ^{-1}(H_i) \rightarrow H_i$
of $\xi$  is smooth and   $\xi ^{-1}(H_i)$  is open in $A$; \item
the union of the  $\xi ^{-1}(H_i)$ generates $A$; \item $A$ is
locally sectionable;  \end{enumerate} \noindent then there is a
unique morphism $\xi ^\prime : A \rightarrow Hol$  of Lie
groupoids such that $\psi \xi ^\prime = \xi $   and $\xi ^\prime h
= i\xi h$ for $h \in \xi^{-1}(H_i)$.
\end{thm}

The groupoid $Hol$ is called the {\it holonomy groupoid}
$Hol(s,{\cal U}_s)$ of the   Lie  local subgroupoid $s$ and atlas
${\cal U}_s$.

We outline  the proof of which full details are given in
\cite{Ao-Br}. Some details of part of the construction  are needed
for Proposition \ref{locex}.

\noindent {\bf Outline proof:}\

Let $G=glob(s)$ and let $\Gamma (G)$ be the set of all admissible
local sections of $G$.  Define a product on $\Gamma (G)$ by
\[(ts)x = (t\beta sx)(sx)\] for two  admissible local sections $s$
and $t$. If $s$ is an admissible local  section then write
$s^{-1}$ for the admissible local section $\beta  s{\cal D}(s)
\rightarrow G, \beta sx \mapstochar \rightarrow (sx)^{-1}$.  With
this product $\Gamma (G)$ becomes an inverse semigroup. Let
$\Gamma ^{r}(W)$ be the subset of $\Gamma (G)$ consisting of
admissible local sections which have values in $W$ and are smooth.
Let $\Gamma ^{r}(G, W)$ be the subsemigroup of $\Gamma (G)$
generated by $\Gamma ^{r}(W)$. Then $\Gamma ^{r}(G, W)$ is again
an inverse semigroup. Intuitively, it contains information on the
iteration of local procedures.

Let $J(G)$ be the sheaf of germs of admissible local sections of
$G$.  Thus the elements of $J(G)$ are the equivalence classes of
pairs $(x,s)$  such that $s \in \Gamma (G), x \in {\cal D}(s)$,
and $(x,s)$ is equivalent to $(y,t)$ if and only if $x = y$ and
$s$ and $t$ agree on a neighbourhood of  $x$. The equivalence
class of $(x,s)$ is written $[s]_{x}$. The product  structure on
$\Gamma (G)$ induces a groupoid structure on $J(G)$ with $X$  as
the set of objects, and source and target maps $[s]_{x}
\mapstochar \rightarrow x, [s]_{x} \mapstochar \rightarrow \beta
sx$. Let $J^{r}(G, W)$  be the subsheaf of $J(G)$ of germs of
elements of $\Gamma ^{r} (G,W)$.  Then $J^{r} (G, W)$ is generated
as a subgroupoid of $J(G)$ by the sheaf  $J^{r} (W)$ of germs of
elements of $\Gamma ^{r} (W)$. Thus an element of $J^{r} (G, W)$
is of the form \[[s]_{x} = [s_n]_{x_n} \ldots [s_1]_{x{_1}}\]
where $s = s_n \ldots s_1$ with $[s_i]_{x_{i}} \in J^{r}(W),
x_{i+1}  = \beta s_{i}x_{i} , i = 1,\ldots ,n$ and $x_1 = x \in
{\cal D}(s) $.

Let $\psi : J(G) \rightarrow G$ be the final map defined by $\psi
([s]_x) = s(x) $, where $s$ is an admissible local section.  Then
$\psi (J^{r} (G, W)) = G $. Let $J_0 = J^{r} (W) \cap \ker \psi $.
Then $J_0$ is a normal subgroupoid of $J^{r} (G, W) $; the proof
is  the same as in~\cite{Ao-Br} Lemma 2.2. The holonomy groupoid
$Hol =  Hol(G, W)$ is defined to be the quotient $J^{r} (G, W)/J_0$.
Let $p:J^{r}(G, W) \rightarrow Hol$ be the quotient morphism and
let $p([s]_{x})$ be denoted by $<s>_{x}$.  Since $J_{0} \subseteq
\ker \psi $ there is a surjective morphism  $\phi : Hol \rightarrow
G$ such that $\phi p = \psi $.\par

The topology on the holonomy groupoid $Hol$ such that $Hol$ with
this topology  is a Lie groupoid is constructed as follows.  Let
$s \in \Gamma ^{r}(G, W)$. A partial function $\sigma _s :  W
\rightarrow Hol$ is defined as follows. The domain of $\sigma _s$
is the set of $w \in W$ such that $\beta w \in {\cal D}(s)$. A
smooth admissible local section $f$ through $w$ is chosen and the
value $\sigma _{s}w$ is defined to be \[\sigma _{s}w =
<s>_{\beta w}<f>_{\alpha w} = <sf>_{\alpha w}.\] It is proven
that $\sigma _{s}w$ is independent of the choice of the local
section $f$ and that these $\sigma _{s}$ form a set of charts.
Then the initial topology with respect to the charts $\sigma _{s}$
is imposed on  $Hol$. With this topology $Hol$ becomes a Lie
groupoid. Again the proof is essentially the same as in
Aof-Brown~\cite{Ao-Br}.

We now outline the proof of the universal property.

Let $a \in A$. The aim is to define $\xi'(a) \in Hol$.

Since $\xi^{-1}(W)$ generates $A$ we can write $a=a_n\ldots a_1$
where $\xi(a_i) \in W$ and hence $\xi(a_i) \in H_{i'}$ for some
$i'$. Since $A$ has enough continuous admissible  local  sections,
we can choose continuous admissible local sections  $f_i$   of
$\alpha _A $ through $a_i , i=1,\ldots,n$,  such that  they  are
composable and their images are contained in  $\zeta \io
(H_{i'})$. The smoothness of $\xi$ on $\xi \io (W)$ implies that
$\xi f_i$ is a smooth admissible local section  of $\alpha$
through  $\xi a_i \in H_{i'}$ whose image is contained in
$H_{i'}$. Therefore $\xi f\in \Gamma^c (G,W)$. Hence we can set
$$\xi 'a=\langle \xi f\rangle_{\A a} \in Hol.$$

The major part of the proof is in showing that $\xi'$ is well
defined, smooth, and is the unique such morphism. We refer again
to \cite{Ao-Br}.

\begin{rem} {\em The above construction shows that the
holonomy groupoid  $Hol (G,W)$ depends on the class ${\cal C}^r$
chosen, and so should strictly  be written $Hol^{r}(G,W)$. An
example of this dependence is given  in Aof-Brown~\cite{Ao-Br}.}
\end{rem}

From the construction of the holonomy groupoid we easily obtain
the  following extendability condition.

\begin{prop} The locally Lie groupoid  $(G,W)$
is extendable to a Lie  groupoid structure on  $G$  if and only if
the following condition holds:  \par \noindent {\bf (1)}:  if $x
\in O_{G}$, and $s$  is a product $s_{n} \ldots s_1$  of local
sections about $x$  such that each $s_i$  lies in $\Gamma ^{r}(W)$
and $s(x) = 1_x$, then there is a restriction $s^\prime $ of $s$
to a neighbourhood of $x$  such that $s^\prime $ has image in $W$
and is  smooth, i.e. $s^\prime \in \Gamma ^{r}(W)$.\label{locex}
\end{prop}

\begin{pf} The canonical morphism $\phi : H \rightarrow G$ is an
isomorphism if and only if $\ker \psi \cap J^{r}(W) = \ker \psi $.
This is equivalent to $\ker \psi \subseteq J^{r}(W)$.  We now show
that $\ker \psi \subseteq J^{r}(W)$  if and only if the condition
(1) is satisfied.\par

Suppose $\ker \psi \subseteq J^{r}(W)$. Let $s = s_n \ldots s_1$
be a product of admissible local sections about $x \in O_G$ with
$s_i \in \Gamma^{r}(W)$ and $x \in {\cal D}_s$ such that $s(x) =
1_x$.  Then $[s]_x \in J^{r}(G,W)$ and $\psi ([s]_x) = s(x) =
1_x$.  So $[s]_x \in \ker \psi$, so that $[s]_x \in J^{r}(W)$. So
there is a neighbourhood $U$ of $x$ such  that the restriction
$s\mid U \in \Gamma ^{r}(W)$.\par

Suppose the condition (1) is satisfied. Let  $[s]_x \in  \ker
\psi$. Since $[s]_x \in J^{r}(G,W)$,  then $[s]_x = [s_n]_{x_n}
\ldots  [s_1]_{x_1}$  where $s = s_n \ldots  s_1$ and
$[s_i]_{x_{i}} \in J^r(W)$, $x_{i+1} = \beta s_{i}x_{i}, i=1,
\ldots ,n$  and $x_1 = x \in {\cal D}(s)$.  Since $s(x) = 1_x$,
then by (1), $[s]_x \in J^{r}(W)$.  \end{pf}

In effect, Proposition \ref{locex} states that the
non-extendability of $(G,W)$  arises from the {\it holonomically
non trivial} elements of $J^{r} (G,W) $.  Intuitively, such an
element $h$ is an iteration of local procedures (i.e. of elements
of $J^{r}( W)$) such that $h$ returns to the  starting point (i.e.
$\alpha h = \beta  h$) but $h$ does not return to the starting
value (i.e. $\psi h \ne 1$). \par

The following gives a circumstance in which this extendability
condition  is easily seen to apply.

\begin{cor}[Corollary 4.6 in \cite{Br-Mu1}]
\label{locexcor} {\it Let} $Q$ {\it be} a Lie {\it groupoid  and
let} $p : M \rightarrow Q$ {\it be a morphism of groupoids such
that} $p : O_{M} \rightarrow O_{Q}$ {\it  is the identity. Let}
$W$ {\it be an open subset of} $Q$ {\it such that}
\begin{enumerate}[a)]
\item   $O_{Q} \subseteq W$;
\item  $W = W^{-1}$;
\item   $W$  generates $Q$;
\item $ ( \alpha _{W}, \beta _{W}, W )$ is smoothly locally
sectionable;
\end{enumerate}
\noindent and suppose that $\tilde{\imath} : W \rightarrow M$  is
given such that $p \tilde{\imath} = i : W \rightarrow Q$  is the
inclusion and  $W^\prime = \tilde{\imath}(W)$ generates $M $.

Then $M$  admits a unique structure of Lie groupoid such that
$W^\prime $ is an open subset and  $p : M \rightarrow Q$  is a
morphism of Lie groupoids  mapping $W^\prime $ diffeomorphically
to $W $.
\end{cor}

\begin{pf} It is easy to check that $(M,W^\prime )$ is a locally
Lie groupoid. We prove that condition (1)  in Proposition
\ref{locex} is satisfied (with $(G,W)$ replaced by $(M,W^\prime )
)$.\par

Suppose given the data of (1). Clearly, $ps = ps_n \ldots ps_1$,
and so $ ps $ is smooth, since $G$ is a Lie groupoid. Since $s(x)
= 1_x$, there is a restriction $s^\prime $ of $s$ to a
neighbourhood of $x$ such that $Im(ps) \subseteq W $. Since $p$
maps $W^\prime$ diffeomorphically to $W $,  then $s^\prime $ is
smooth and has image contained in $W $. So (1) holds, and by
Proposition \ref{locex}, the topology on $W^\prime $ is extendable
to make  $M$ a Lie groupoid. \end{pf}

\begin{remark}  {\em It may seem unnecessary to construct the
holonomy groupoid in order to verify extendability  under
condition (1) of Proposition \ref{locex}. However the construction
of the smooth  structure on $M$ in the last corollary, and the
proof that this yields  a Lie groupoid, would have to follow more
or less the steps  given in Aof and Brown~\cite{Ao-Br} as sketched
above. Thus it is more  sensible to rely on the general result. As
Corollary \ref{locexcor} shows,  the utility of (1) is that it is
a checkable condition, both positively  or negatively, and so
gives clear proofs of the non-existence or existence  of
non-trivial holonomy. }
\end{remark}

Putting  everything together gives immediately our main theorem on
monodromy.
\begin{thm}(Strong Monodromy Principle)\label{smp}
Let $s$ be a strictly regular Lie local subgroupoid of a Lie
groupoid $Q$, and let ${\cal U}_s=\{(U_i, H_i): i\in I\}$ be a
strictly regular atlas  for  $s$. Let  $W({\cal
U}_s)=\bigcup_{i\in I}H_i$. Then there is a Lie groupoid
$M=Mon(s,{\cal U}_s)$ and morphism $p: M\to glob(s)$ which is the
identity on objects with the following properties:
\begin{enumerate}[a)] \item The injections $H_i\to glob(s)$ lift to
injections $\eta_i: H_i\to M$ such that $W^\prime =\bigcup_{i\in
I}\eta_i(H_i)$ is an open submanifold of $M$.
\item  $W^\prime$ generates $M$,
\item If $K$ is a Lie groupoid and $f =\{ f_i:H_i\to K, i\in
I\}$ is a smooth local morphism, then there is a unique smooth
morphism $M(f):M\to K$ extending the $f_i, i\in I.$
\end{enumerate}
\end{thm}
\begin{pf}
Starting with $s$ we form the locally Lie groupoid $(glob(s),
W({\cal U}_s)$ and then its holonomy groupoid $Hol(glob(s),
W({\cal U}_s))$. Regarding $W({\cal U}_s))$ as contained in $Hol$
we can form the monodromy groupoid $M=M(W({\cal U}_s))$ with its
projection to $Hol$. By Corollary \ref{locexcor} (with $Q=Hol$) M
obtains the structure of Lie groupoid.

Conditions a) and b) are immediate from this construction of the
monodromy groupoid.

In c), the existence of $M(f)$ follows from the weak monodromy
principle. To prove that $M(f)$  is smooth it is enough, by local
sectionability, to prove  it is smooth    at the identities  of $M
$.   This  follows  since $p \colon M\rightarrow G$    maps
$\tilde{\imath}(W)$ diffeomorphically to  $W $.
\end{pf}

\begin{rem}{\em
We have now formed from a strictly regular Lie local subgroupoid $s$
of the Lie groupoid $G$ a smooth morphism of Lie groupoids
\[  \xi :Mon(s, {\cal U}_s)\to Hol(s, {\cal U}_s) \]
which is the identity on objects so that the latter holonomy
groupoid is a quotient of the monodromy groupoid. It also follows
from \cite[Proposition 2.3]{Br-Mu1} that this morphism is a
covering map on each of the stars of these groupoids.

Extra conditions are needed to ensure that $\xi$ is a universal
covering map on stars -- see   \cite[Theorem 4.2]{Br-Mu1}. This
requires further investigation, for example we may need to shrink
$W$ to satisfy the required condition.

This also illustrates that Pradines' theorems in \cite{Pr1} are
stated in terms of germs. Again, the elaboration of this needs
further work.
 }\end{rem}

\begin{rem}{\em  The above results also include the notion of a
Lie local equivalence relation, and a strong monodromy principle
for these. We note also that the Lie groupoids we obtain are not
\'{e}tale groupoids. This is one of the distinctions between the
direction of this work and that of Kock and Moerdijk
\cite{Ko-Mo1,Ko-Mo2}. It would be interesting to investigate the
relation further, particularly with regard to the monodromy
principle.

A further point is that a local equivalence relation determines a
topos of sheaves of a particular type known as an \`etendue
\cite{Ko-Mo2}. What type of topos is determined by a local
subgroupoid?
 }
\end{rem}

{}

\end{document}